\documentclass[12pt]{article}
\usepackage{a4wide}
\newcommand{\nn}{\nonumber \\}
\begin{document}
\title{A note on the $q$-derivative operator}
\author{J. Koekoek \and R. Koekoek}
\date{ }
\maketitle

\begin{abstract}
The following statement is proved. If the $q$-derivative operator ${\cal D}_q$ is defined by
$${\cal D}_qf(x):=\left\{\begin{array}{ll}\displaystyle\frac{f(x)-f(qx)}{(1-q)x},& x\ne 0\\  \\
f'(0),& x=0\end{array}\right.$$
for functions $f$ which are differentiable at $x=0$, then we have for every
positive integer $n$
$$\left({\cal D}_q^nf\right)(0)=\lim\limits_{x\rightarrow 0}{\cal D}_q^nf(x)=
\frac{f^{(n)}(0)}{n!}\frac{(q;q)_n}{(1-q)^n}$$
for every function $f$ whose $n$th derivative at $x=0$ exists.

We give a proof in both the real variable and the complex variable case.
\end{abstract}

\vspace{5mm}

Usually the $q$-derivative operator ${\cal D}_q$ is defined by
$${\cal D}_qf(x)=\frac{f(x)-f(qx)}{(1-q)x}$$
where $q$ is fixed and $0<q<1$. However, this definition is not valid for $x=0$.
On the other hand, the hypothesis $0<q<1$ can be weakened. So we define
\begin{equation}
\label{def}
{\cal D}_qf(x):=\left\{\begin{array}{ll}\displaystyle\frac{f(x)-f(qx)}{(1-q)x},& x\ne 0\\  \\
f'(0),& x=0\end{array}\right.
\end{equation}
for functions $f$ which are differentiable at $x=0$, where $q$ is a fixed real
number not equal to 1.

Note that the existence of $f'(0)$ implies that the domain of $f$ contains
the point zero in its interior.

Further we define ${\cal D}_q^nf:={\cal D}_q\left({\cal D}_q^{n-1}f\right)$ for $n=1,2,3,\ldots$,
where ${\cal D}_q^0$ denotes the identity operator.

In this paper we prove the following theorem.

\vspace{5mm}

{\bf Theorem.} {\sl Let $n$ be a positive integer and let $f$ be a function
for which $f^{(n)}(0)$ exists. Then we have
$$\left({\cal D}_q^nf\right)(0)=\lim\limits_{x\rightarrow 0}{\cal D}_q^nf(x)=
\frac{f^{(n)}(0)}{n!}\frac{(q;q)_n}{(1-q)^n}.$$}

\vspace{5mm}

{\bf The real variable case.}

\vspace{5mm}

Let $f$ be a function of a real variable $x$ and let the domain of $f$ contain
the interval $(-\rho,\rho)$ for some $\rho>0$.

We use the following version of l'Hospital's rule
(see for instance section 84 in \cite{Frank})~:

\vspace{5mm}

{\sl Let $F$ and $G$ be functions whose $n$th ($n\ge 1$) derivatives
at $x=0$ exist. Suppose that
$$F^{(m)}(0)=G^{(m)}(0)=0,\;m=0,1,2,\ldots,n-1\;\mbox{ and }\;G^{(n)}(0)\ne 0.$$
Then we have
$$\lim\limits_{x\rightarrow 0}\frac{F(x)}{G(x)}=\frac{F^{(n)}(0)}{G^{(n)}(0)}.$$}

\vspace{5mm}

First of all we take $q=0$. Then we have
$${\cal D}_0f(x):=\frac{f(x)-f(0)}{x},\;0<|x|<\rho$$
and $\left({\cal D}_0f\right)(0):=f'(0)$. Further we have
\begin{eqnarray}
\label{A}
\left({\cal D}_0^nf\right)(0)&:=&\left({\cal D}_0^{n-1}f\right)'(0)\nn
&=&\lim_{x\rightarrow 0}\frac{{\cal D}_0^{n-1}f(x)-\left({\cal D}_0^{n-1}f\right)(0)}{x}
=\lim_{x\rightarrow 0}{\cal D}_0^nf(x),\;n=1,2,3,\ldots.
\end{eqnarray}
Now we will show that
\begin{equation}
\label{3'}
x^n{\cal D}_0^nf(x)=f(x)-\sum_{k=0}^{n-1}\frac{f^{(k)}(0)}{k!}x^k,\;0<|x|<\rho
\end{equation}
and
\begin{equation}
\label{5'}
\left({\cal D}_0^nf\right)(0)=\frac{f^{(n)}(0)}{n!}
\end{equation}
for each $n\in\{1,2,3,\ldots\}$ for which $f^{(n)}(0)$ exists.

We use induction on $n$. For $n=1$ we have (\ref{3'}) and (\ref{5'}) by
definition. So we assume that $f^{(n+1)}(0)$ exists and that (\ref{3'}) and
(\ref{5'}) are valid for some $n\in\{1,2,3,\ldots\}$. Then we find
\begin{eqnarray}
x^{n+1}{\cal D}_0^{n+1}f(x)&:=&x^n\left[{\cal D}_0^nf(x)-\left({\cal D}_0^nf\right)(0)\right]\nn
&=&f(x)-\sum_{k=0}^{n-1}\frac{f^{(k)}(0)}{k!}x^k-\frac{f^{(n)}(0)}{n!}x^n\nn
&=&f(x)-\sum_{k=0}^n\frac{f^{(k)}(0)}{k!}x^k,\;0<|x|<\rho.\nonumber
\end{eqnarray}
Further we have by using (\ref{A}) and l'Hospital's rule above
\begin{eqnarray}
\left({\cal D}_0^{n+1}f\right)(0)&=&\lim_{x\rightarrow 0}{\cal D}_0^{n+1}f(x)\nn
&=&\lim_{x\rightarrow 0}\frac{f(x)-\sum\limits_{k=0}^n\frac{f^{(k)}(0)}{k!}x^k}{x^{n+1}}
=\frac{f^{(n+1)}(0)}{(n+1)!}.\nonumber
\end{eqnarray}
This proves (\ref{3'}) and (\ref{5'}) for each $n\in\{1,2,3,\ldots\}$ for
which $f^{(n)}(0)$ exists.

Now a combination of (\ref{A}) and (\ref{5'}) completes the proof of the
theorem if $q=0$.

\vspace{5mm}

For $q=-1$ the definition (\ref{def}) leads to
$$2x{\cal D}_{-1}f(x)=f(x)-f(-x),\;0<|x|<\rho$$
and $\left({\cal D}_{-1}f\right)(0)=f'(0)$. So we have by using l'Hospital's
rule
$$\lim_{x\rightarrow 0}{\cal D}_{-1}f(x)
=\lim_{x\rightarrow 0}\frac{f(x)-f(-x)}{2x}
=f'(0)=\left({\cal D}_{-1}f\right)(0)$$
if $f'(0)$ exists. This proves the theorem for $n=1$.

Further we have
$$4x^2{\cal D}_{-1}^2f(x)=f(x)-f(-x)+f(-x)-f(x)=0,\;0<|x|<\rho.$$
And by using l'Hospital's rule we obtain
\begin{eqnarray}
\left({\cal D}_{-1}^2f\right)(0)&:=&\left({\cal D}_{-1}f\right)'(0)\nn
&=&\lim_{x\rightarrow 0}\frac{{\cal D}_{-1}f(x)-\left({\cal D}_{-1}f\right)(0)}{x}
=\lim_{x\rightarrow 0}\frac{f(x)-f(-x)-2xf'(0)}{2x^2}=0\nonumber
\end{eqnarray}
if $f''(0)$ exists. This implies that
$${\cal D}_{-1}^2f(x)=0,\;|x|<\rho.$$
Hence
$${\cal D}_{-1}^nf(x)=0,\;|x|<\rho,\;n=2,3,4,\ldots.$$
So we have for $n\ge 2$ : if $f^{(n)}(0)$ exists, then $f''(0)$ exists and
${\cal D}_{-1}^nf(x)=0$ for $|x|<\rho$. This proves the theorem for $q=-1$.

\vspace{5mm}

In the sequel we assume that $q\ne 0$ and $q\ne -1$.

\vspace{5mm}

We show first that for every positive integer $n$
\begin{equation}
\label{first}
\lim\limits_{x\rightarrow 0}{\cal D}_q^nf(x)=\frac{f^{(n)}(0)}{n!}\frac{(q;q)_n}{(1-q)^n}
\end{equation}
if $f^{(n)}(0)$ exists.

The definition (\ref{def}) of the $q$-derivative operator ${\cal D}_q$ leads to
$$(1-q)x{\cal D}_qf(x)=f(x)-f(qx),\;0<|x|<\min(\rho,\frac{\rho}{|q|}).$$
If we multiply by $(1-q)x$ and replace $f$ by ${\cal D}_qf$ we obtain
$$(1-q)^2x^2{\cal D}_q^2f(x)=f(x)-(1+q^{-1})f(qx)+q^{-1}f(q^2x),\;0<|x|<\min(\rho,\frac{\rho}{|q|^2}).$$
Now we use induction on $n$ to see that we have for $n\in\{1,2,3,\ldots\}$~:
\begin{equation}
\label{dqn}
(1-q)^nx^n{\cal D}_q^nf(x)=\sum_{k=0}^n(-1)^kq^{-k(n-1)+{k\choose 2}}\left[n\atop k\right]_q
f(q^kx),\;0<|x|<\min(\rho,\frac{\rho}{|q|^n})
\end{equation}
where the $q$-binomial coefficient $\left[n \atop k\right]_q$ is defined by
$$\left[n \atop k\right]_q:=\frac{(q;q)_n}{(q;q)_k(q;q)_{n-k}}.$$

Let $n$ be a positive integer and assume that $f^{(n)}(0)$ exists.

Now we define
$$F(x):=\sum\limits_{k=0}^n(-1)^kq^{-k(n-1)+{k \choose 2}}
\left[n \atop k\right]_qf(q^kx)$$
and
$$G(x):=(1-q)^nx^n.$$
Then we have
$$G^{(m)}(0)=0,\;m=0,1,2,\ldots,n-1\;\mbox{ and }\;G^{(n)}(0)=(1-q)^nn!\ne 0$$
and
$$F^{(m)}(0)=f^{(m)}(0)\sum\limits_{k=0}^n(-1)^kq^{k(m-n+1)+{k \choose 2}}
\left[n \atop k\right]_q,\;m=0,1,2,\ldots,n.$$
Now we use (see for instance \cite{GasRah})
$$(a;q)_n=\sum_{k=0}^nq^{{k \choose 2}}\left[n\atop k\right]_q(-a)^k$$
for $q$ real with $|q|\ne 1$ and $a$ an arbitrary complex number, to obtain
\begin{equation}
\label{sumdelta}
\sum_{k=0}^n(-1)^kq^{k(m-n+1)+{k\choose 2}}\left[n\atop k\right]_q
=(q^{m-n+1};q)_n,\;m=0,1,2,\ldots.
\end{equation}
Hence
$$\sum_{k=0}^n(-1)^kq^{k(m-n+1)+{k\choose 2}}\left[n\atop k\right]_q
=(q;q)_n\delta_{mn},\;m=0,1,2,\ldots,n.$$
So we have
$$F^{(m)}(0)=0,\;m=0,1,2,\ldots,n-1\;\mbox{ and }\;F^{(n)}(0)=(q;q)_nf^{(n)}(0).$$
Applying l'Hospital's rule above we obtain
$$\lim\limits_{x\rightarrow 0}{\cal D}_q^nf(x)=\lim\limits_{x\rightarrow 0}\frac
{\sum\limits_{k=0}^n(-1)^kq^{-k(n-1)+{k\choose 2}}\left[n\atop k\right]_qf(q^kx)}{(1-q)^nx^n}
=\frac{f^{(n)}(0)}{n!}\frac{(q;q)_n}{(1-q)^n}.$$
This proves (\ref{first}).

\vspace{5mm}

Now we prove that for every positive integer $n$
\begin{equation}
\label{second}
\left({\cal D}_q^nf\right)(0)=\frac{f^{(n)}(0)}{n!}\frac{(q;q)_n}{(1-q)^n}
\end{equation}
if $f^{(n)}(0)$ exists.

We use induction on $n$.

For $n=1$ we have (\ref{second}) by definition.

Assume that $f^{(n+1)}(0)$ exists and (\ref{second}) holds for some $n\in\{1,2,3,\ldots\}$.
Now we use the definition (\ref{def}) and (\ref{dqn}) to find
\begin{eqnarray*}
\left({\cal D}_q^{n+1}f\right)(0)&:=&\left({\cal D}_q^nf\right)'(0)=
\lim\limits_{x\rightarrow 0}\frac{{\cal D}_q^nf(x)-\left({\cal D}_q^nf\right)(0)}{x}\\
&=&\lim\limits_{x\rightarrow 0}\frac{\sum\limits_{k=0}^n(-1)^kq^{-k(n-1)+{k\choose 2}}
\left[n\atop k\right]_qf(q^kx)-\frac{f^{(n)}(0)}{n!}(q;q)_nx^n}
{(1-q)^nx^{n+1}}=:L.
\end{eqnarray*}
Let
$$F(x):=\sum\limits_{k=0}^n(-1)^kq^{-k(n-1)+{k\choose 2}}\left[n\atop k\right]_qf(q^kx)
-\frac{f^{(n)}(0)}{n!}(q;q)_nx^n$$
and
$$G(x):=(1-q)^nx^{n+1}.$$
Then we have by using (\ref{sumdelta})
$$F^{(m)}(0)=G^{(m)}(0)=0,\;m=0,1,2,\ldots,n$$
and
$$F^{(n+1)}(0)=(q^2;q)_nf^{(n+1)}(0)\;\mbox{ and }\;G^{(n+1)}(0)=(1-q)^n(n+1)!\ne 0.$$
So we have by using l'Hospital's rule again
$$L=\frac{F^{(n+1)}(0)}{G^{(n+1)}(0)}=\frac{(q^2;q)_nf^{(n+1)}(0)}{(1-q)^n(n+1)!}
=\frac{f^{(n+1)}(0)}{(n+1)!}\frac{(q;q)_{n+1}}{(1-q)^{n+1}}.$$
This proves (\ref{second}) for every $n\in\{1,2,3,\ldots\}$ for which $f^{(n)}(0)$
exists.

\vspace{5mm}

{\bf The complex variable case.}

\vspace{5mm}

Now let $f$ be a function of a complex variable $z$.

If $f'(0)$ exists, then we have by definition $\left({\cal D}_qf\right)(0)=f'(0)$.
Further we have for $q\ne 0$
\begin{eqnarray*}
\lim\limits_{z\rightarrow 0}{\cal D}_qf(z)&=&\lim\limits_{z\rightarrow 0}
\frac{f(z)-f(qz)}{(1-q)z}\\
&=&\lim\limits_{z\rightarrow 0}\left(\frac{f(z)-f(0)}{(1-q)z}-
q\frac{f(qz)-f(0)}{(1-q)qz}\right)
=\frac{f'(0)}{1-q}-q\frac{f'(0)}{1-q}=f'(0).
\end{eqnarray*}

For $q=0$ we simply obtain
$$\lim_{z\rightarrow 0}{\cal D}_0f(z)=\lim_{z\rightarrow 0}\frac{f(z)-f(0)}{z}=f'(0).$$

If $f''(0)$ exists, then $f'$ exists in a neighbourhood of $z=0$
which means that $f$ is analytic in a neighbourhood of $z=0$.
This implies that $f^{(m)}(0)$ exists for $m=0,1,2,\ldots$ and we have
$$f(z)=\sum_{m=0}^{\infty}\frac{f^{(m)}(0)}{m!}z^m,\;|z|<\rho\;\mbox{ for some }\;\rho>0.$$
Applying the definition (\ref{def}) we find
$${\cal D}_qf(z)=\sum_{m=1}^{\infty}\frac{f^{(m)}(0)}{m!}\frac{1-q^m}{1-q}z^{m-1},\;|z|<\min(\rho,\frac{\rho}{|q|}),$$
where $\min(\rho,\frac{\rho}{|q|})$ should be replaced by $\rho$ in the case that $q=0$.

Hence ${\cal D}_qf(z)$ is analytic in a neighbourhood of $z=0$. Now we use
induction to find for $n\in\{1,2,3,\ldots\}~:$
$${\cal D}_q^nf(z)=\sum_{m=n}^{\infty}\frac{f^{(m)}(0)}{m!}\frac{(q^{m-n+1};q)_n}{(1-q)^n}z^{m-n},\;|z|<\min(\rho,\frac{\rho}{|q|^n}),$$
where $\min(\rho,\frac{\rho}{|q|^n})$ should be replaced by $\rho$ if $q=0$.

Hence ${\cal D}_q^nf(z)$ is analytic in a neighbourhood of $z=0$. So we easily obtain
$$\lim\limits_{z\rightarrow 0}{\cal D}_q^nf(z)=\left({\cal D}_q^nf\right)(0)=\frac{f^{(n)}(0)}{n!}\frac{(q;q)_n}{(1-q)^n}.$$
This completes the proof of the theorem in the complex variable case.

\vspace{5mm}

{\bf Remark 1.} Since we defined the $q$-derivative operator for functions $f$
which are only assumed to be differentiable at $x=0$, the well-known limit
$$\lim\limits_{q\rightarrow 1}{\cal D}_qf(x)=f'(x)$$
no longer holds in general, but only for those $x$ where $f$ is differentiable.

\vspace{5mm}

{\bf Remark 2.} In the complex variable case the definition (\ref{def}) can
be extended to complex values of $q$ and then the theorem still holds for
complex values of $q\ne 1$.

\vspace{5mm}

\small

\noindent
J. Koekoek\\Menelaoslaan 4\\5631 LN Eindhoven\\The Netherlands

\vspace{5mm}

\noindent
R. Koekoek\\Delft University of Technology\\Faculty of Technical Mathematics
and Informatics\\Mekelweg 4\\2628 CD Delft\\The Netherlands

\end{document}